\makeatletter
\def\@maketitle{%
\defaultfont\normalsize
\let\@makefnmark\relax \let\@thefnmark\relax \ifx\@empty\@subjclass\else
\@footnotetext{1991 {\it Mathematics Subject 
Classification}.\enspace
\@subjclass.}\fi
\ifx\@empty\@keywords\else
\@footnotetext{{\it Key words and phrases.}\enspace \@keywords.}\fi
\ifx\@empty\@thanks\else
\@footnotetext{\@thanks}\fi
\topskip66\p@ 
\vtop{\centering{\baselineskip14\p@\bf
\expandafter{\@title}\@@par}%
\global\dimen@i\prevdepth}%
\prevdepth\dimen@i
\ifx\@empty\@authors
\else
\baselineskip32\p@
\vtop{\@andify{ AND }\@authors
\centering{{\@authors}\@@par}%
\global\dimen@i\prevdepth}\relax
\prevdepth\dimen@i
\fi
\ifx\@empty\@dedicatory
\else
\baselineskip18\p@
\vtop{\centering{\small\it\@dedicatory\@@par}%
\global\dimen@i\prevdepth}\prevdepth\dimen@i \fi
\ifx\@empty\@date\else
\baselineskip24\p@
\vtop{\centering\@date\@@par
\global\dimen@i\prevdepth}\prevdepth\dimen@i \fi
\normalsize
\dimen@32\p@ \advance\dimen@-\baselineskip \vskip\dimen@\@plus14\p@
} 
\makeatother

\documentstyle[amscd,amssymb,verbatim,amsthm,12pt]{amsart}

\theoremstyle{plain}
\newtheorem{Thm}{Theorem}
\newtheorem*{Thm*}{Theorem}
\newtheorem{Cor}[Thm]{Corollary}
\newtheorem{Lem}[Thm]{Lemma}
\newtheorem{Prop}[Thm]{Proposition}

\theoremstyle{remark}

\newtheorem{Rem}[Thm]{Remark}

\newcommand{\ra}{\rightarrow}
\newcommand{\Spec}{\operatorname{Spec}}

\renewcommand{\O}{{\Cal O}}

\def\uint#1#2{{U\kern -8pt\int}_{\kern -5pt {#1}}^{\kern 2pt {#2}} }
\def\lint#1#2{{L\kern -7.5pt \int}_{\kern -5pt {#1}}^{\kern 2pt {#2}} }
\def\morph#1{\overset{#1}{\ra}}

\def\dsp#1{$\displaystyle{#1}$}
\def\Cal#1{{\cal #1}}
\def\BBB#1{{\Bbb #1}}

\def\set#1{{\{{#1}\}}}

\theoremstyle{remark}

\numberwithin{equation}{section}
\renewcommand{\rm}{\normalshape}

\newcommand{\A}{{\Cal A}}

\renewcommand{\O}{{\Cal O}}

\newcommand{\D}{{\Cal D}}

\newcommand{\Pic}{\operatorname{Pic}}
\newcommand{\Jac}{\operatorname{Jac}}

\newcommand{\Ber}{\operatorname{Ber}}

\newcommand{\C}{\BBB C}

\def\dv{\text{Div}}

\def\comp#1#2{{#1}\kern-2pt-\kern-2pt{#2}}

\def\C{{\BBB C}}

\begin{document}

\title[Serre duality, Abel's theorem, Jacobi inversion]{Serre duality, Abel's theorem, and Jacobi inversion for supercurves over a thick superpoint}

\author{Mitchell J. Rothstein}
\address{Department of Mathematics, University of Georgia, Athens, GA 30602, USA}
\email{rothstei@@math.uga.edu}

\author{Jeffrey M. Rabin}
\address{Department of Mathematics, UCSD, La Jolla, CA 92093, USA}
\email{jrabin@@math.ucsd.edu}

\begin{abstract}The principal  aim of this paper is to extend Abel's theorem to the setting of complex supermanifolds of dimension $1|q$ over a finite-dimensional local supercommutative $\BBB  C$-algebra.  The theorem is proved by  establishing a compatibility of Serre duality for the supercurve with Poincar\'e duality on the reduced curve.   We include an elementary algebraic proof of the requisite form of Serre duality, closely based on the account of the reduced case given by  Serre  in {\em Algebraic groups and class fields},  combined with an  invariance result for the topology on the dual  of the space of r\' epartitions.  Our Abel map,  taking  Cartier divisors of degree zero to the dual of the space of sections of the Berezinian sheaf, modulo  periods,  is defined via Penkov's characterization of the Berezinian sheaf as the cohomology of the de Rham complex of the sheaf $\D$ of differential operators,  as a right module over itself.  We discuss the Jacobi inversion problem for the Abel map and give an example demonstrating that if $n$ is an integer sufficiently large that the generic divisor of degree $n$ is linearly equivalent to an effective divisor, this need not be the case for {all} divisors of degree $n$.
\end{abstract}
\maketitle

\section{Introduction}

 In the classical theory of Riemann surfaces, a fundamental role is played by the Abel map, which links the algebraic theory of projective curves with the transcendental theory of Riemann surfaces \cite{GH}.
Abel's theorem states that two divisors of  degree zero are linearly equivalent if and only if they have the same image under the Abel map.  In this paper we prove that this statement remains valid for supercurves of dimension $1|q$ over a thickened point, by which we mean $\Spec(B)$, where 
$B$ is a a finite-dimensional local supercommutative $\BBB C$-algebra.  Part of the task is to define the Abel map in this setting.  This was done for Weil divisors with $q=1$ in \cite{BR}.
Here we give a definition for arbitrary $q$ using Cartier divisors.  To construct the target of the Abel map, we use the characterization of the Berezinian sheaf, $\Ber$,  as the cohomology of the de Rham complex of the sheaf $\D$ of differential operators on the structure sheaf, $\O$ \cite{P}.
The period map, equation \eqref{eq:oldper} below,  maps $H_1(X,\BBB Z)\to H^0(X,\Ber)^\circ$. \footnote{
If $R$ is a $\BBB Z_2$-graded ring, all $R$-modules shall be tacitly assumed to be   $\BBB Z_2$-graded.
We define $Hom$ in the category of $\BBB Z_2$-graded $R$-modules in such a way that  $Hom_R(M,N)$ consists of  parity-preserving $R$-module homomorphisms.     Thus ``maps" are parity-preserving by default.
By definition, automorphisms preserve parity.
If $R$ is supercommutative,  we also have the {\em internal } hom functor,  adjoint to the tensor product,  denoted   $\underline{Hom}_R(M,N)$.  Then
${Hom}_R(M,N)$ is the even part of $\underline{Hom}_R(M,N)$.  We define the {\em dual} of $M$ to be the internal hom,   $\underline{Hom}_R(M,R)$,  and denote it by  $M^*$. 
Its even part ${Hom}_R(M,R)$ will be denoted $M^\circ$.}
Defining $\Pic^0(X)$ as the group of divisors of  degree zero modulo linear equivalence,
and $\Jac(X)$ as the quotient $H^0(X,\Ber)^\circ/H_1(X,\BBB Z)$, Abel's theorem then says that the Abel map imbeds $\Pic^0(X)$ in $\Jac(X)$.  The infinitesimal version of this statement,
that $H^1(X,\O)$ imbeds in $H^0(X,\Ber)^*$, is a corollary of Serre duality, theorem \ref{thm:dual},
which in the present setting says that  $H^0(X,\Ber)=H^1(X,\O)^*$.\footnote{This raises the question of whether $H^1(X,\O)$ is reflexive,  i.e., isomorphic to its double dual.  This is guaranteed if $B$ is Gorenstein, and in particular if $B$ is a Grassmann algebra \cite{D,E}.We do not have an example of a supercurve for which $H^1(X,\O)$ is not reflexive.}  The classical Jacobi inversion theorem,  asserting that every divisor of degree  equal to the genus of $X$ is linearly equivalent to an effective divisor, has no immediate analogue for arbitrary supercurves, which depend topologically on both the genus $g(X)$ of the reduced space  and the Chern class $c(X)$ of the vector bundle associated to $\O_X$. One possible extension of the assertion is that if
$n$ is an integer sufficiently large that the generic divisor of degree $n$ is linearly equivalent to an effective divisor, then the same holds for all divisors of degree $n$.  The compactness argument that yields this statement in the classical case is not available here,  and we give an example to show that it is false in the $(1|1)$-dimensional case.

\section{Serre duality}

Serre duality was established for complex supermanifolds in \cite{HW},  for projective supervarieties over a field in \cite{OP} and \cite{P}, and for projective supervarieties over a Grassmann algebra in \cite{BR}.  It seems to be  ``known" in the more general setting of derived categories and morphisms of superschemes, though we are not aware of a reference.  We therefore include  an elementary proof of 
the version we need here, following the proof of the classical case found in \cite{S}.

To begin, let $ B$ be a finite-dimensional local supercommutative $ \C$-algebra. 
Fix a positive integer $q$.  By the term ``smooth $1|q$ dimensional supercurve over $B$" (or simply ``supercurve" if $ B$ and $q$ are understood) we shall mean a pair
$X=(X_0,\O)$, where $X_0$ is a smooth projective curve over $\BBB C$, and $\Cal O$ is a sheaf of supercommutative $B$-algebras in the Zariski topology,  such that there is a cover of $X_0$ by  open sets   $\set{U_{\alpha}}$ satisfying for every index $\alpha$, 
\begin{equation}\label{eq:lt}
\Cal O|_{U_{\alpha}}\simeq 
  B[\theta_1,\ldots,\theta_q]\otimes_{ \C}\O_0|_{U_{\alpha}}\end{equation}
in the category of sheaves of local supercommutative $ B$-algebras.  Here $\O_0$ denotes the algebraic structure sheaf of $X_0$ and  
$B[\theta_1,\ldots, \theta_q]$ is the Grassmann algebra over $B$ on {\it q} generators.

Let $X=(X_0,\O)$ be a supercurve. The sheaf of meromorphic functions on $X$ is by definition the sheaf of fractions  $f/g$,  where $f\in\O$,  $g\in\O$, and $g$ is even and not nilpotent.
 The {\em function ring} of $X$,  $B(X)$, is the ring of global meromorphic functions.   As in the non-super case,  for every point $P\in X_0$,  $B(X)$ is isomorphic to the ring of fractions of the local ring $\O_P$.

Keeping $X_0$ fixed throughout the discussion,  let $X[B,q]$ denote the trivial family $(X_0,\O[B,q])$, where globally $$\O[B,q]=
  B[\theta_1,\ldots,\theta_q]\otimes_{ \C}\O_0 .
$$
Let 
$$\Lambda[B,q]=
  B[\theta_1,\ldots,\theta_q]\otimes_{ \C}\C(X_0).
$$
From the local triviality  \eqref{eq:lt} it follows that (non-canonically)  $B(X)$ is isomorphic to 
$\Lambda[B,q]$.

Denote by $n(R)$ the nilpotent ideal of an arbitrary supercommutative  ring $R$, and by
$Aut^+(R)$ the kernel of the natural map $Aut(R)\to Aut(R/n(R))$.
 Then we have, for every point $P\in X_0$,
$$Aut(\O[B,q]_P)\subset Aut^+(\Lambda[B,q]).$$
Denoting by $\A ut(\O[B,q])$ the automorphism sheaf of $\O[B,q]$,  we therefore have an inclusion of sheaves
$$\A ut(\O[B,q])\to Aut^+(\Lambda[B,q]).$$

Let $\D[B,q]$ denote the sheaf of linear differential operators on $\O[B,q]$.

\begin{Lem}\label{lem:subD}
$\A ut(\O[B,q])\subset \D[B,q].$
\end{Lem}

\begin{pf}
Let $\tau\in Aut^+(\Lambda[B,q])$. Then
$$\tau(\theta_i)=\alpha_i+\sum_j A_{ij}\theta_j+\cdots$$
where $\alpha_i$ and $A_{ij}$ belong to $ B\otimes \BBB C(X)$ and the ellipsis denotes terms of higher degree in $\theta_j$. The $ B\otimes \BBB C(X)$-linear map sending $\theta_i$ to $\alpha_i+\sum A_{ij}\theta_j$ determines an automorphism of $\Lambda[B,q]$, and is a differential operator.  After composing with the inverse of this automorphism,  we may assume that $id-\tau$ maps  $\Lambda[B,q]$ to the ideal generated by
the nilpotents in $B$ and the square of the nilpotents in $\Lambda[B,q]$.  Letting $Z$ denote $id-\tau$,  $Z$ satisfies $Z(fg)=fZ(g)+Z(f)g-Z(f)Z(g)$. It follows by induction that $Z$ is a nilpotent differential operator. 
\end{pf}

 For any sheaf of groups $\Cal S$,  let
$\Pi^0\Cal S\subset\Pi_{P\in X_0}\Cal S_P$ denote the set of elements $\eta$ such that $\eta_P$ is the identity element for all but finitely many $P$.  Let
 $\gamma\in\Pi^0  Aut^+(\Lambda[B,q])$.  One obtains a supercurve $X^{\gamma}=(X_0,\O^{\gamma})$ by taking $\O^{\gamma}$ to be the subsheaf of $\Lambda[B,q]$ such that for all $P$, $\O^{\gamma}_P=\gamma_P(\O[B,q]_P)$. 

By local triviality, \eqref{eq:lt}, one has
\begin{Prop}\label{prop:lt}
All supercurves are of the form $X^\gamma$ for some element $\gamma$.
\end{Prop}

Let $\Lambda[B,q]^\times $ denote the group of even units. For all $\xi\in\Pi^0\Lambda[B,q]^\times$ we get a rank-one locally free sheaf of $\O^ {\gamma}$-modules as follows: Let
$\O^{\gamma}(\xi)$ denote the subsheaf of $\Lambda[B,q]$ such that for all $P$, $\O^{\gamma}(\xi)_P=\gamma_P(\xi_P\O[B,q]_P)$.   (As usual, $\O^{\gamma}(\xi)$ depends only on the divisor class of $\xi$,  but this divisor class will depend on $\gamma$.)
Once again, every rank-one locally free sheaf on $X$ is of this form.
  
A r\' epartition  on $\O[B, q]$ is 
 a map $r:X_0\to  \Lambda[B,q]$ such that $r_P\in\O[B,q]_P$ for all but finitely many $P$ (cf. Serre \cite{S}). 
  Let $R[B,q]$ denote the set of all r\'epartitions. Regard $\Lambda[B,q]$ as   a subring of $R[B,q]$, identifying $\Lambda[B,q]$ with constant functions.
 Define the subset $R(\gamma,\xi)\subset R[B,q]$ as the set of functions $r$ such that for all points $P$, $r_P\in \O^{\gamma}(\xi)_P$.
Then as in \cite{S}, Prop. II.3, 
\begin{equation}\label{eq:h1}
H^1(X_0,\O^{\gamma}(\xi))\simeq R[B,q]/(R(\gamma,\xi)+\Lambda[B,q]).
\end{equation}

For fixed $\gamma$, let $R[B,q]$ be given the topology such that the spaces  $R(\gamma,\xi)$ for all $\xi$ form a neighborhood base at $\set 0$.  
Then $H^1(X_0,\O^{\gamma}(\xi))^*$ is the annihilator of $R(\gamma,\xi)$ in the topological dual of  $R[B,q]\,/ \Lambda[B,q]$.

\begin{Prop}\label{prop:top} The topology on
$R[B,q]$ is independent of $\gamma$.

\end{Prop}

\begin{pf}  Let $\sigma:X_0\to\A ut^+(\Lambda[B,q])$ be another finitely supported function.   By lemma \ref{lem:subD}, $\gamma_P$ and $\sigma_P$ are meromorphic differential operators. It follows that if
 $t_P\in\O_0(P)$ is a local parameter at $P$, there exists an integer $m_P$ such that $ \sigma_P t_P^{m_P}\gamma_P^{-1}$ is regular at $P$ as a differential operator.  Then
 $R(\gamma,\xi)\subset R(\sigma,\tau)$, where $\tau_P=t_P^{-m_P}$ for $P$ in the support of $\gamma$ or $\sigma$, and $\tau_P=1$ elsewhere. 
\end{pf}

\begin{Thm}\label{thm:continuous} Let $\omega_0$ be a nonzero  meromorphic one-form  on $X_0$.  Then   each continuous element of $(R[B,q]\,/ \Lambda[B,q])^*$ is of the form 
\begin{equation}\label{eq:pair}
g\mapsto\sum_P res_P(\omega_0 \, \partial_{\theta_1} \cdots \partial_{\theta_q}(fg))
\end{equation}
 for a unique $f\in\Lambda[B,q]$. In particular, every element of $H^1(X_0,\O^{\gamma}(\xi))^*$ is of this form.
\end{Thm}

\begin{pf}
Note first that
$$R[B,q]\,/ \Lambda[B,q]\simeq B[\theta _1,\ldots,\theta_q]\otimes_{\BBB C} (R(X_0)/\BBB C(X_0))$$  
where $R(X_0)$ is the space of r\'epartitions on the reduced space.

The topological $\BBB C$-linear dual of $R(X_0)/\BBB C(X_0)$ is the space of meromorphic one-forms on $X_0$, via the residue pairing \cite{S}.  The $B$-linear dual of $B[\theta _1,\ldots,\theta_q]$ is itself, via  the pairing
$$ f\cdot g= \partial_{\theta_1} \cdots \partial_{\theta_q}(fg).$$
The theorem follows from the definition of the tensor product.
\end{pf}

 Given a supercurve $X$,  let $\Ber_X$ denote the Berezinian of its cotangent sheaf.
%

\begin{Thm}[Serre Duality for supercurves]\label{thm:dual}
Let $X$ be a supercurve and let $\Cal L$ be a rank-one locally free sheaf on $X$. 
Then the formula
\begin{equation}\label{eq:SD}g\mapsto\sum_P res_P(dz \, \partial_{\theta_1} \cdots \partial_{\theta_q}(fg))
\end{equation}
 defines a pairing of
$g \in H^1(X_0,\Cal L)$ with $f \in H^0(X_0,\Ber_X \otimes \Cal L^{-1})$, with respect to which
$$H^1(X_0,\Cal L)^*=H^0(X_0,\Ber_X\otimes \Cal L^{-1}).$$ 
\end{Thm}
\begin{pf}
The sense in which \eqref{eq:SD}
 defines a pairing of
$H^1(X_0,\Cal L)$ with $H^0(X_0,\Ber_X \otimes \Cal L^{-1})$ will  be made clear in the course of the proof. 
We may assume that $\O_X=\O^{\gamma}$ and $\Cal L=\O^{\gamma}(\xi)$.
Then $H^1(X_0,\O^{\gamma}(\xi))^*$ is the annihilator of $R(\gamma,\xi)$.
Let $f\in\Lambda[B,q]$.  Then formula \eqref{eq:pair} defines an element of $H^1(X_0,\O^{\gamma}(\xi))^*$ if and only if $f$ satisfies a set of local conditions.  If $\xi=1$,  the  conditions are that for all $P$, and all $h\in\Lambda[B,q]$,  if $\gamma_P(h)\in\O_P$, then 
$$res_P(\omega_0 \, \partial_{\theta_1} \cdots \partial_{\theta_q}(fh))=0.$$
Having chosen $\omega_0$ arbitrarily, we may redefine $f$ so that $\omega_0=dz$ for some local parameter $z\in(\O_0)_P$.

The change of variables formula for Berezin integration is as follows:
Let
$z_1,\ldots,z_{p},\theta_1,\ldots,\theta_{q}$ and $w_1,\ldots,w_{p},\eta_1,\ldots,\eta_{q}$
be two coordinate systems near a point $P$
on a $p|q$-dimensional supermanifold.  Let $dz\, \partial_\theta$ denote the $\Omega^p$-valued differential operator $dz^1\wedge \cdots \wedge dz^p\otimes\partial_{\theta_1} \cdots \partial_{\theta_q}$.
Then
\begin{equation}\label{eq:cov}
dz \, \partial_\theta \, Ber\binom{\partial w\ \partial\eta}{\partial z\ \partial\theta}=dw \, \partial_\eta+\epsilon 
\end{equation}
where $$
\epsilon=d\circ L$$ for some $\Omega^{p-1}$-valued differential operator $L$.   (This is the statement that Berezin integration is well-defined modulo boundary terms, \cite{R}.)

Let $w=\gamma_P^{-1}(z)$, $\eta_i=\gamma_P^{-1}(\theta_i)$.
Then
\begin{align}\label{eq:resop2}
res_P \circ dz \, \partial_\theta (f)=
res_P \circ dw \, \partial_\eta \Ber\binom{\partial z\ \partial\theta}{\partial w\ \partial\eta}
(f).
\end{align}
This shows that formula \eqref{eq:pair} defines an element of
$H^1(X_0,\Cal L)^*$ if  and only if the meromorphic section $\omega_0 \, \partial_{\theta_1} \cdots\partial_{\theta_q}$ of $\Ber_{X[B,q]}$ is a holomorphic section of $\Ber_{X^{\gamma}}$. The rest of the theorem follows
from theorem \ref{thm:continuous}.

\end{pf}
                                                                                              
%
%
%
%
%

Theorem \ref{thm:dual} does {not} guarantee that 
$H^1(X,\Cal L)$ is the dual of  $H^0(X,\Ber_X \otimes {\Cal L}^{-1})$ without further conditions on $B$.  (See footnote 2). It is known that the cohomology groups are finitely generated, so we do have

\begin{Cor}\label{cor:inject}
The  pairing \eqref{eq:pair} gives an injection 
$$0\to H^1(X,\Cal L)\rightarrow  H^0(X,\Ber_X\otimes\Cal L^{-1})^*.$$
\end{Cor}

\begin{Rem} Classically, Serre duality has the following corollary for a Riemann surface $X_0$
(see, for example, Corollary 4.4 in \cite{G}). A differential principal part $\frak p$  extends to a meromorphic differential on $X_0$ if and only if
$\sum_{P \in X_0} res_P ({\frak p}) = 0$.  Taking $\Cal  L=\Ber$ in corollary \ref{cor:inject},  we get the natural generalization of this result.
%
\begin{Prop} There is a meromorphic section of $\Ber$ on $X$ having a given principal part 
$\frak p$ if and only if 
$$ \sum_{P \in X_0} res_P (\partial_{\theta_1} \cdots \partial_{\theta_q} ({\frak p} g)) = 0 $$
for every global holomorphic function $g \in H^0(X,\cal O)$.  \end{Prop}
Note that, in general, global sections $g$ of $\cal O$ need not be constant,  and $H^0(X,\cal O)$ need not be freely generated.\end{Rem}

\section{Abel's Theorem}

In this section we work  in the complex topology:  $\O$ now stands for $\O_{hol}$.  Besides making the Poincar\'e lemma available \cite{DM,V}, this will give us an interpretation of the residue on $X$ as the pairing of  a section of $\Ber$ defined on an annulus with the fundamental class of the annulus.  This pairing  also gives rise to the period map.

\subsection{Definition of the Abel map}

We begin by reviewing Penkov's characterization of $\Ber$, \cite{P}.
Let $\Cal D_k$ denote the sheaf of $B$-linear differential operators from $\O$ to the sheaf of $k$-forms $\Omega^k$.  Write $\Cal D$ for $\Cal D_0$.  One has the de Rham complex
$$\cdots \to
\Cal D_k\to\Cal D_{k+1}\to \cdots$$
given by  $L\mapsto d\circ L$, where $d$ is the exterior derivative.
Note that this is a complex of $\O$-modules under right multiplication.
Penkov observes that
 the first (or more generally $p^{th}$ for a supermanifold of dimension $p|q$) cohomology sheaf of this complex is rank-one locally free, with basis $dz \, \partial_{\theta_1} \cdots \partial_{\theta_q}$,  where $(z,\theta_j)$ are local coordinates, the other cohomology sheaves being $0$.  This a strengthened form of the change of variables formula.

Thus,  letting $\Cal D_{1,cl}\subset\Cal D_1$ denote the subsheaf
consisting of differential operators $L:\O\to\Omega^1$  such that for all $f\in\O$,  $d\circ L(f)=0$, one has an exact sequence
$$0\to\Cal D\morph{d\circ}\Cal D_{1,cl}\morph{\pi}
\Ber\to 0.$$

Let $U\subset  X_0$ be an annulus,  and let  $[U]\in H_1(U,\BBB Z)$ be its fundamental class.  
Let $\omega$ be a section of $\Ber$ on $U$.  We wish to recover the residue as a pairing of $\omega$ with $[U]$, which we will then denote  by $\oint_{[U]}\omega$.  There are (at least) two methods to define it.
\bigskip

Method 1:  First define the residue of a closed one-form.
By the Poincar\'e lemma,  we have exactness of
$$0\to B\to \O\to \Omega^1_{cl}\to 0.$$
Define $\oint_{[U]}$ on  $H^0(U,\Omega^1_{cl}
)$ to be the connecting homomorphism $H^0(U,\Omega^1_{cl})\to H^1(U,B)=B$.  Then, because $U$ is Stein, we may lift $\omega$ to a section $L\in H^0(U,\D_{1,cl})$ and define
\begin{equation}\label{eq:res1}
\oint_{[U]}\omega=\oint_{[U]} L(1).
\end{equation}
\bigskip

Method 2:  
  let $\Cal D^\sharp\subset
\Cal D_{1,cl}$ denote the kernel of the map
$$L\mapsto L(1).$$

\begin{Lem}\label{lem:onto}
The map $\Cal D^\sharp\morph{\pi}\Ber$ is surjective.
\end{Lem}

\begin{pf}
Let $L$ be a section of $\Cal D_{1,cl}$ on an open set $U$ and let $\omega=\pi(L)$.  $L(1)$ is closed, and therefore by shrinking $U$ if necessary we can assume there exists a section $f$ of $\O$ such that $L(1)=df$.  Regarding $f$ as a section of $\Cal D$,  we then have
$$\omega=\pi(L-d\circ f). $$
This proves the claim.
\end{pf}
Let $\Cal D^\flat\subset\Cal D$ denote the subsheaf of $\Cal D$ that annihilates the constant sheaf $B$.  (That is,  $\Cal D^\flat$ is  the left ideal generated by vector fields.)  By lemma \ref{lem:onto} we have a short exact sequence
\begin{equation}\label{eq:ses}0\to B\oplus \Cal D^\flat\morph {\large{d}\overset{\circ}{}}\Cal D^\sharp\morph{\pi}
\Ber\to 0. \end{equation}
We may then define $\oint_{[U]}\omega$ to be image of $\omega$ under the connecting homomorphism
of \eqref{eq:ses} on $U$, $$H^0(U,\Ber)\to H^1(U,B)\oplus H^1(U,\D^{\flat}),
$$ followed by 
projection onto $H^1(U,B)$.
 
\begin{Lem}Methods 1 and 2 give the same result.  Moreover,  if $\omega$ is meromorphic,  then the residue as defined in terms of Laurent series satisfies
$$res_P(\omega)=\oint_{[U]}\omega$$
where $U$ is a sufficiently small deleted neighborhood of $P$.\end{Lem}
\noindent In view of this, we will sometimes write the residue as $\oint_P \omega$.

Similarly, the 
{\em period  map}
\begin{equation}\label{eq:oldper} H^0(X,\Ber)\morph{per} H^1(X,B)
 \end{equation}
is defined to be the connecting homomorphism
of \eqref{eq:ses} on $X$, $$H^0(X,\Ber)\to H^1(X,B)\oplus H^1(X,\D^{\flat})
$$ followed by 
projection onto $H^1(X,B)$.
It can be shown that in the case $q=1$, this is the period map in \cite{BR}.

Note:  By Serre duality, $H^0(X,\Ber)$ is a dual module, and is therefore \cite{D} naturally isomorphic to its double dual. So there will be no information lost if we dualize the period map and continue to call it $per$.

We now have a diagram

\begin{equation}\label{eq:commutes}
\begin{CD}
H^1(X,B)   @>{i}>>   H^1(X,\O)   \\
@VVV                                    @VV{res}V       \\
H_1(X,B)   @>{per}>>   H^0(X,\Ber)^*
\end{CD}
\end{equation}

\noindent
where the left arrow is Poincar\'e duality, $res$ denotes Serre duality, and the top arrow comes from the inclusion $B\to\O$.

\begin{Lem}\label{lem:commutes}
Diagram \eqref{eq:commutes} commutes.
\end{Lem}

\begin{pf}
Let $\set{A_1,\dots,A_g,B_1,\dots,B_g}\subset H_1(X,\BBB Z)$ be a standard homology basis, with dual basis $\set{A^1,\dots,A^g,B^1,\dots,B^g}\subset H^1(X,\BBB Z)$.  It is enough to check commutativity on a basis for $H^1(X,B)$; as a representative case we choose $A^1 \in H^1(X,B)$ and verify that $per(B_1)=res(i(A^1))$. Let $\omega\in H^0(X,\Ber)$.  Represent $A_1$ and $B_1$ in the standard way as embedded circles intersecting at one point.  Let $P$ be a point disjoint from $A_1$ and $B_1$.  Let $V_0$ be a small disk containing $P$ and let $V_1=X_0-P$.  With respect to the open cover $\set{V_0,V_1}$, the image of $A^1$ in $H^1(X,\O)$ is represented by a section
 $g\in H^0(V_0-P,\O)$. Let $U_1\subset X_0$ be an annulus containing $B_1$ and let $U_0=X_0-B_1$. 
Let $L_0\in H^0(V_0,\Cal D_{1,cl})$ and $L_1\in H^0(U_1,\Cal D_{1,cl})$ be representatives of $\omega$ on $V_0$ and $U_1$ respectively.
 We must show  that
$$\oint_{[U_1]}L_1(1)=\oint_{[V_0-P]}L_0(g).$$
  The intersection $U_1\cap U_0$ is a pair of disjoint annuli,  $W_{\pm}$.  The \u Cech cocycle representing $A^1$ with respect to the open cover $\set{U_0,U_1}$ is the function on
 $U_1\cap U_0$ that equals $1$ on $W_+$ and $0$ on $W_-$.  Then
 the statement that the image of $A^1$ in $H^1(X,\O)$ is represented by
 $g$ is expressed in terms of \u Cech cohomology by the following two statements:
 
 1. $g$ extends to $X_0-(\set{P}\cup B_1)$.
 
 2. There is a section $f\in H^0(U_1,\O)$ such that $g-f$ is $1$ on $W_+$ and $0$ on $W_-$.
 
 Because $V_1$ is Stein, the section $L_1$ representing $\omega$ on the annulus $U_1$ can be chosen such that it is defined on all of $V_1$.  Then
 \begin{equation}\label{eq:homology}
\oint_{[V_0-P]} L_0(g) =\oint_{[V_0-P]} L_1(g).
\end{equation}
The boundary of $U_0-V_0$ consists of three circles,  namely $[V_0-P]$ employed in equation \eqref{eq:homology} and $\partial U_1$ consisting of two circles homologous to $[W_\pm]$.  Thus
 \begin{equation}\label{eq:transfer}
\oint_{[V_0-P]} L_1(g) = \oint_{[W_+]}L_1(g-f+f)-\oint_{[W_-]}L_1(g-f+f).
\end{equation}
But $f$ is defined on all of $U_1$, so
$$\oint_{[W_+]}L_1(f)-\oint_{[W_-]}L_1(f)=0. $$
Therefore
$$
\oint_{[V_0-P]} L_1(g) =
\oint_{[W_+]} L_1(g-f)-\oint_{[W_-]} L_1(g-f)=\oint_{[W_+]}L_1(1)
$$
which is what we needed to show, since $[W_+]$ is homologous to $[U_1]$.
\end{pf}

\begin{Rem} In \cite{BR} the map $res \circ i$ was called $rep$. The formula for $rep$ in terms of a canonical homology basis stated without proof in Lemma 2.9.1 of \cite{BR} follows from the above computation.
\end{Rem}

Define the sheaf of Cartier divisors $\dv X$ by the exact sequence

\begin{equation}\label{eq:div}
0\to\O_X^\times\to\Cal K_X^\times\to \dv X\to 0,
\end{equation}
where $\Cal K^\times$ is the sheaf of invertible even meromorphic sections of $\O$
and $\O^\times=\Cal K^\times\cap\O$.

If $P$ is a point in $X_0$ and
 $f\in({\Cal K_X^\times})_P$,  one has the quantity $\oint_P\frac{df}f$, which depends only on the class of $f$ in $({\dv X})_P$.  Just as in the non-super case, one has 

\begin{Lem}\label{lem:integer}
\dsp{\oint_P \frac{df}f} is an integer.
\hfill\qed\break
\end{Lem}

To define the Abel map,  let $\xi\in H^0(X,\dv X)$ such that the {\it degree} 
$\sum_P \oint_P\xi=0$.  Let $\omega\in H^0(X,\Ber)$.  Fix a connected simply connected open set $U$ containing the support of $\xi$.  Choose $L\in H^0(U,\Cal D^\sharp)$ representing
 $\omega$,  and choose $f\in H^0(U,\Cal K^\times)$ representing $\xi$.
 For each point $Q\in U$ not belonging to the support of $\xi$ there is a germ $g\in\Cal O_Q$ such that $e^g=f$ in a  neighborhood of $Q$.  Since $L$ annihilates constants,  $L(g)$ is independent of which logarithm of $f$ is chosen.  Thus we obtain a
 closed one-form defined on $U-\{$singular points of $f\}$, and we 
 may unambiguously write this one-form as
 $L(\log f)$.
 
Let us examine the quantity 
 $$\rho=\sum_{P\in U}\oint_PL(\log f) \in B$$
 with regard to the choices made.
At each point $P$,  $f$ may be altered by multiplying by $e^h$ for some $h$ in the even part of $\Cal O_P$.  This changes $L(\log f)$ to $L(\log f)+L(h)$,  which does not alter $\rho$.  The section $L$ may be altered by adding $d\circ (c+M)$ for some constant $c\in B$ and section $M\in H^0(U,\Cal D^\flat)$.    In a neighborhood of $P$ we may choose vector fields $Y_i$ and differential operators $M_i$ such that $M=\sum M_i Y_i$. Writing 
$f=e^g$ as before,  we have $$d\circ M(g)=d\bigg(\sum M_i\bigg(\frac{Y_i(f)}f\bigg)\bigg)$$ which is annihilated by $\oint_P$.
Therefore
$$\sum_{P\in U}\oint_P(L+d\circ(c+M))(\log f)=\sum_{P\in U}\oint_PL(\log f)+ c\sum_{P\in U}\oint_Pdf/f$$
and the second term on the right-hand side vanishes by hypothesis.  We therefore have a well-defined element of $B$,  independent of the representatives of $\xi$ and $\omega$,  but depending however on the open set $U$.  Let $U_1=U$,  and let $U_2$ be another connected simply connected open set containing the support of $\xi$.  We then must consider
$$\sum_{P\in U_1}\oint_PL_1(\log f_1)-\sum_{P\in U_2}\oint_PL_2(\log f_2).$$
The independence of this quantity on the representatives of $\xi$ remains valid.  However, if we let $W_1,\dots, W_n$ denote the connected components of $U_1\cap U_2$, then 
$L_1-L_2=d\circ (c+M)$,  where $c$ takes a constant value $c_i$ on each $W_i$.  We have integers $n_i=\sum_{P\in W_i}\oint_P\xi$,  summing to $0$.
Then $$
\sum_{P\in U_1}\oint_P L_1(\log f_1)-\sum_{P\in U_2}\oint_P L_2(\log f_2)=\sum_i c_i n_i$$
which is the pairing of the class of $\omega$ in $H^1(X_0,B)$ with a homology class in $H_1(X_0,\BBB Z)$.  

We therefore have a map, the {\em Abel map},
$$H^0(X,\dv X)_0\morph{\bf Abel} Hom(H^0(X,\Ber),B)/H_1(X,\BBB Z) = 
H^0(X,\Ber)^\circ/H_1(X,\BBB Z)$$
mapping the group of divisors of total degree zero to the Jacobian of $X$.

\subsection{Abel's Theorem}

Let $U\subset X_0$ be a disk, let $H^0(U,\dv X)_0$ denote the sections of $\dv X$ over $U$ having total degree zero,  and let
$${\bf ab}_U: H^0(U,\dv X)_0\to H^0(X,\Ber)^\circ$$
denote the map  described in the previous section.  We have the following diagram:

\begin{equation}\label{eq:Abel}
\begin{CD}  H^1(X,\O)_{ev} @>\text{exp}>>H^1(X,\O^\times) & \\
@V{res}VV                            @A{\ell}AA&\\
   H^0(X,\Ber)^\circ@<{{\bf ab}_U}<<H^0(U,\dv X)_0&\subset H^0(X,\dv X)_0
\end{CD}
\end{equation}
where $\ell$ is the connecting homomorphism for the sequence \eqref{eq:div}.

\begin{Lem}\label{lem:factor}   Let $\gamma\in H^1(X,\O)_{ev}$ and $\xi\in H^0(U,{\rm \dv} X)_0$.  If  $res(\gamma)={\bf ab}_U(\xi)$,  then $exp(\gamma)=\ell(\xi)$.

%
%

\end{Lem}
\begin{pf}
Let $\gamma\in H^1(X,\O)_{ev}$ and $\xi\in H^0(U,{\rm \dv} X)_0$ such that $res(\gamma)={\bf ab}_U(\xi)$.
Represent $\xi$ by a meromorphic section $f\in H^0(U,\Cal K^\times)$.  Represent $\gamma$ by a meromorphic function $g$ defined on $U$, with a pole at one point $Q$ not belonging to the  support of $\xi$.  Let $D_i$,  $i=1,2$ be disjoint disks contained in $U$, such that $Q\in D_1$ and $supp(\xi)\subset D_2$.  From the fact that $\xi$ has degree zero it follows that there is a branch of $\log f$ defined on $U-D_2$.

Let $\omega\in H^0(X,\Ber)$.  Represent $\omega$ by $L\in H^0(U,\D^\sharp)$.
We are given that for all $\omega$,
$$\sum_{P\in U}\oint_P L(\log f)=res_Q(L(g)).$$
Let $U'\subset U$ be a slightly smaller disk, such that  $U-U'$ is an annulus and $D_1\cup D_2\subset U'$. Then the previous equation can be rewritten as 
$$\oint_{[U-U']} L(g-\log f)=0.$$
This holds for all $\omega$,  so by 
Cor. \ref{cor:inject}  the cohomology class in $H^1(X,\O)$ defined by the section $g-\log f \in H^0(U-U',\Cal O)$ with respect to the open cover $\set{X-U',U}$ is equal to zero.  Thus there exist sections $h^-\in H^0(X-U',\Cal O)$ and $h^+\in H^0(U,\Cal O)$ such that on $U-U'$,
$$\log(f)-g=h^--h^+.$$
Then we obtain an invertible section of $\O$, $$\phi\in H^0(X-(\set{Q}\cup(\text{support of\ }\xi)),\O^\times)$$
 by patching together $e^{h^-}$ on $X-U'$ and $e^{h^+}e^{-g} f$ on $U$.  Letting $\O_g$ denote the line bundle with transition function $e^g\in H^0(D_1-Q,\O^\times)$,  and letting $\O(\xi)$ denote the line bundle associated to the divisor $\xi$,  we see that $\phi$ is a trivialization of $\O_g^{-1}\otimes\O(\xi)$,  which is what we needed to show.\end{pf} 

 \begin{Cor}[{\bf Abel's Theorem}]\label{cor:abel}
Let $\xi\in H^0(X,{\rm \dv} X)_0$,  and let $\O(\xi)$ be the associated line bundle.  Then $\O(\xi)$ is trivial, i.e, $\xi$ is the divisor of a globally defined section $f\in H^0(X,\Cal K^\times)=B(X)$, if and only if ${\bf Abel}(\xi)=0$.
\end{Cor}

\begin{pf}
Assume ${\bf Abel}(\xi)=0$,  and let $U$ be a disk containing the support of $\xi$.  Then there exists $c\in H_1(X,\BBB Z)$  such  that ${\bf ab}_U(\xi)=per(c)$.  Let $\gamma\in H^1(X,\BBB Z)\subset H^1(X,\O)$ be the image of $c$ under Poincar\'e duality. 
By lemma \ref{lem:commutes},  $res(\gamma)={\bf ab}_U(\xi)$. Thus $\O(\xi)$ is trivial, by lemma \ref{lem:factor}.  The converse is proved by the classical argument, \cite{GH}: If $f$ is a global meromorphic function,  
then the Abel image of the divisor class of $a+bf$ depends homogeneously on $[a,b] \in \BBB P^1$ and is therefore constant.\end{pf}

\section{Jacobi Inversion}

A divisor $\xi$ is {\em effective} if at each point $P$, it can be represented locally by a function  
$f\in \Cal K^\times\cap\O_P$.  
Given an arbitrary divisor $\xi$,  the existence of an effective divisor $\xi'$ linearly equivalent to 
$\xi$ is equivalent to the existence of a non-nilpotent section of $\O(\xi)$. A perturbative argument with respect to the nilpotent ideal
%
%
 shows that
  for $n$ sufficiently large,  every divisor  of degree $n$ is linearly equivalent to an effective divisor.   Let $n(X)$ denote the least such $n$.  One may also consider the least $n$ such that a generic
 divisor  of degree $n$ is linearly equivalent to an effective divisor.   Denote this number by  $n_{gen}(X)$. In the classical case, $q=0$, the Jacobi inversion theorem asserts that $n(X)=g$, the genus of $X_0$,  \cite{GH}.  The proof  first establishes that $n_{gen}(X)=g$  and then uses compactness to  prove that $n_{gen}(X)=n(X)$.  The following example shows that $n(X)$ may be strictly larger than $n_{gen}(X)$ when $q>0$.

 Let $\Cal L$  be  a generic line bundle of degree $0$ on $X_0$.  Let $Y$ denote the $1|1$ dimensional supercurve over $\Spec\BBB C$ with odd part $\Cal L$ and let $X=Y\times_{\Spec \BBB C}\Spec\BBB C[\beta]$, where $\beta$ is an  odd parameter.
Given a class $c\in H^1(X_0, \Cal L)$,  let $\Cal F_c$ denote the line bundle on $X$ with transition data $1+\beta c\in H^1(X,\O^*)$. 
 Let $\Cal J$ be a line bundle on $X_0$, of degree $g$,  such that 
 $h^0(\Cal J)=1$ and $h^0(\Cal J\otimes \Cal L)>1$.  Note that by Riemann-Roch,
 $h^1(\Cal J\otimes \Cal L)>0$. Our supermanifold is split,  so we can pull back  $\Cal J$ to $X$.  That is to say, on $X$ we have the line bundle $\pi^*(\Cal J)=\Cal J\otimes_{\O_{X_0}}\O_X$, where $\pi:X\to X_0$ is the splitting. Multiplication by $\beta$ gives 
the short exact sequence 
 \begin{equation}\label{eq:1}
 0\to\Cal L\otimes_{\O_{X_0}}\Cal J\morph{\beta}(\pi^*(\Cal J)\otimes_{\O_X}\Cal F_c)_{even}\to \Cal J\to 0\end{equation}
Let $\phi$ be a nonzero section of $\Cal J$.  Then the image of $\phi$ under the connecting homomorphism of \eqref{eq:1} is $\phi c\in H^1(X_0,\Cal L\otimes\Cal J)$.
The map
%
\begin{equation}\label{eq:3}
H^1(X_0,\Cal L)\to H^1(X_0,\Cal L\otimes\Cal J)\end{equation}
sending $c$ to $\phi c$ is surjective,  and in particular it is not the zero map given our choice of $\Cal L $ and $\Cal J$.  If we choose $c$ such that $\phi c\ne 0$,  then the connecting homomorphism of \eqref{eq:1} is injective, and therefore all sections of $\pi^*(\Cal J)\otimes_{\O_X}\Cal F_c$ are nilpotent.  On the other hand,  let  $\Cal M$ be a generic line bundle on $X$ of degree  $g$,  and let $\Cal M_0$ denote its restriction to $X_0$.  Again one has a short exact sequence
$$  0\to\Cal L\otimes_{\O_{X_0}}\Cal M_0\morph{\beta}\Cal M_{even}\to \Cal M_0\to 0$$
For generic $\Cal M$, $H^1(X_0,\Cal L\otimes_{\O_{X_0}}\Cal M_0)=0$. 
Thus we have $n(X)>g$ and $n_{gen}(X)=g$.

\bibliographystyle{alpha}

\end{document}